\documentclass{amsart}
\usepackage[T1]{fontenc}

\theoremstyle{definition}

\theoremstyle{remark}

\numberwithin{equation}{section}

\begin{document}
\title{On the development of differential geometry in Estonia}
\author{Viktor Abramov}
\address{Institute of Pure Mathematics, University of Tartu, Liivi 2, %
Tartu 50409, Estonia}%
\email{abramov@ut.ee}%
\thanks{This work was supported by the grant of Estonian Science Foundation
No. 4515}%
\begin{abstract}%
We give a brief review of a research made in the field of
differential geometry in Estonia in the period from the beginning
of the 19th century to the present time. The biographic data of
mathematicians who made a valuable contribution to the development
of differential geometry in Estonia in mentioned period are
presented. The material of the first section covers the period
from the beginning of the 19th century to the middle of the 20th
century and it can be considered as a brief historical sketch of
the development of differential geometry in Estonia in this
period. The next sections give an idea of the modern trends of
development of differential geometry in Estonia.
\end{abstract}%
\maketitle%
The history of research in the field of differential geometry in
Estonia is inseparably linked with the history of the University
of Tartu. Therefore we begin by reminding some basic facts from
the history of the University of Tartu. The University of Tartu
was founded by the King Gustav II Adolf of Sweden in 1632. We omit
the turbulent period from the foundation of the university to the
end of 18th century when the university sometimes was in Tartu,
sometimes was forced to move to other towns of Estonia such as
Tallinn and P\"arnu as a result of agreements between belligerent
powers. In 1802 the university was reopened in Tartu as a
provincial Baltic university depending upon the local Knighthoods
- it was titled Kaiserliche Universit\"at zu Dorpat (also
Imperatorskij Derptskij Universitet). A first important landmark
in the history of differential geometric research at the
University of Tartu goes back to this period, when J. Martin
Bartels (1769-1836) became a professor of mathematics at the
University of Tartu (1821). Johann Martin Christian Bartels was
born in Braunschweig. He studied at the University of Helmstedt
and then at the G\"ottingen University. He took his doctoral
degree at the University of Jena with a thesis in the field of
variational calculus in 1803. It should be mentioned that Bartels
was the first teacher of Gauss in Braunschweig. Before he was
invited to occupy a professorship of mathematics at the University
of Tartu, Bartels was a professor of mathematics at the University
of Kasan (Russia) for twelve years from 1808 to 1820, where one of
his students was Nicolai Ivanovitch Lobachevsky future professor
of mathematics and rector of the University of Kasan and one of
the founders of non-Euclidean geometry. Bartels contributed to the
theory of space curves by creating the method now known as the
method of moving orthonormal frame. Given a space curve one can
define the local trihedron at a point $P$ of this curve, which
consists of three orthogonal unit vectors ${\bf t}, {\bf n}, {\bf
b}$, where ${\bf t}$ is the unit tangent vector, ${\bf n}$ is the
principal normal vector and ${\bf b}$ is the binormal vector. The
triple $\{{\bf t}, {\bf n}, {\bf b}\}$ bears the name of Frenet
frame at a point of a curve. Bartels studied the rate of change of
the trihedron $\{{\bf t}, {\bf n}, {\bf b}\}$, when a point $P$
begins to move to a near point $Q$ along a curve and he was the
first to derive the equations expressing the derivatives of the
vectors ${\bf t}, {\bf n}, {\bf b}$  in terms of the vectors
themselves now known as Frenet-Serret formulae. The formulae
obtained by Bartels were published by his disciple C.E. Senff in
1831, which means that they appeared 17 years earlier than the
equations published by Frenet and 22 years before Serret published
his equations \cite{Lumiste0}. It should be noted that Bartels
used the components of the vectors ${\bf t}, {\bf n}, {\bf b}$ not
the vectors themselves because the notion of a vector was actually
developed later.

The professorship of applied mathematics was opened in 1843 and
Ferdinand Minding (1806--1885) was invited to occupy this position.
Ferdinand Minding was born in Kalisz (Poland) and shortly after
that his family moved to Hirschberg (now Jelenia Gora in Poland).
He studied classical philology and philosophy at Halle and Berlin
universities from 1824 to 1828. Working as a secondary school
teacher he completed his thesis on approximate calculation of
double integrals and successfully defended it at the University of
Halle in 1829. Ferdinand Minding made a valuable contribution to
the theory of surfaces. He defined the geodesic curvature of a
curve and proved that this curvature is constant along the
shortest curve encircling the given area on a surface. In the
series of papers published from 1838 to 1849 Minding laid the
foundations of the theory of surface bending. In 1864 he was
elected to St. Petersburg Academy of Sciences as an
associate-member and as a honorary member in 1879.

Minding lectured on the theory of curves and surfaces and among
the mathematics students of the University of Tartu attending his
lectures was Karl Peterson (1828--1881). He was born in Riga and
studied mathematics at the University of Tartu from 1847 to 1852.
In 1853 Peterson completed his thesis "On the bending of surfaces"
and having defended it obtained the candidate degree. Upon
graduation, he failed to get a position at the University of
Tartu and moved to Moscow, where he served as a mathematics
instructor at the Petropavlov specialized school.

In spite of an appreciation given by Minding to the Peterson'
thesis "On the bending of surfaces" it was not published until
1952, when appeared the Russian translation of the thesis made by
L. Depman. In Peterson' thesis we find two equations which can be
written in the modern notations as follows
\begin{eqnarray}
\frac{\partial \Delta}{\partial%
     v}-\frac{\partial\Delta'}{\partial u}+\Gamma^{2}_{22}\Delta%
       -2\;\Gamma^{2}_{12}\Delta'+\Gamma^{2}_{11}\Delta''&=&0,\label{firstI}\\%
\frac{\partial \Delta''}{\partial%
     u}-\frac{\partial\Delta'}{\partial v}+\Gamma^{1}_{22}\Delta%
       -2\;\Gamma^{1}_{12}\Delta'+\Gamma^{1}_{11}\Delta''&=&0,\label{secondII}%
\end{eqnarray}
where%
$$
\Delta=\frac{L}{\sqrt{EG-F^{2}}},\quad%
   \Delta'=\frac{M}{\sqrt{EG-F^{2}}},\quad%
     \Delta''=\frac{N}{\sqrt{EG-F^{2}}},
$$
and $E,G,F$ are the coefficients of the first fundamental form
$g$ of a surface;\ \ $L,N,M$ are the coefficients of the second
fundamental
form $h$ of a surface and $\Gamma^{i}_{jk}$ are Christoffel symbols. %
The equations (\ref{firstI}), (\ref{secondII}) play an essential
role in the theory of surfaces. We remind that Gaussian curvature
$K$ of a surface can be written in the form
\begin{equation}
K=\frac{LM-N^{2}}{EG-F^{2}}\label{third}.
\end{equation}
Substituting the Gaussian curvature $K$ in the above formula by
its expression in the terms of the coefficients of the first
fundamental form $g$ we obtain a relation between the coefficients
of the first and second fundamental forms. It turns out that this
relation and relations (\ref{firstI}), (\ref{secondII}) determine
a surface up to congruence. By other words it can be proved that
given two quadratic forms $g, h$, where $g$ is positively
definite, with coefficients satisfying the equations
(\ref{firstI}), (\ref{secondII}), (\ref{third}) there exists a
surface whose first fundamental form is $g$, the second is $h$ and
a surface is determined up to congruence. Shortly after
publication of the Russian translation of the Peterson' thesis it
was generally recognized that Peterson was the first who obtained
the fundamental equations of the theory of
surfaces (\ref{firstI}), (\ref{secondII}) and anticipated the
fundamental theorem of surface theory \cite{Phillips}.%

At the same time that Minding was studying the theory of surfaces,
the subject was also occupying the attention of another scientist
from the University of Tartu, Thomas Clausen (1801--1885). Thomas
Clausen was born in Northern Jutland and having accepted an
invitation to occupy the position of an astronomer-observer, he
came to Tartu in 1842. Clausen interest towards  the theory of
surfaces was inspired by the paper of C.G.J. Jacobi, who tried to
generalize the Gauss theorem about the sum of interior angles of a
geodesic triangle. Casting doubt on the correctness of the results
published by C.G.J. Jacobi, Clausen elaborated a new proof for the
Gauss theorem. Inspired by an another Jacobi paper, where an
integration of differential equation determining the geodesic
lines of an ellipsoid was reduced to quadratures, Clausen showed
that the same reduction can be made in the case of any
non-developed second order surface. Clausen also studied the lunes
of Hippocrates and he showed how two new types of squarable lunes
with proportions of circular arcs 5:1 and 5:3 can be found. It
should be also mentioned that Clausen found a new way of
determining the lemniscate and this work is related to the field
of geometrical constructions.

The successors of Minding on the chair of applied mathematics at
the University of Tartu, who did research in the field of
geometry, were Otto Staude (1857--1928) and Adolf Kneser
(1862--1937). Staude was the first who began to construct the
second order surfaces with the help of the tautened threads. He
also studied the geodesic curvature of a line on a surface and
the sign of the torsion of a curve. Adolf Kneser studied the
algebraic lines by means of synthetic methods and he proved that
if a plane curve has no other singularities except double
tangents and only one double point, then it has only one double
tangent.

Friedrich Schur (1856--1932) took up the post of the professor of
pure mathematics at the University of Tartu in 1888, succeeding on
this post to Peter Helmling (1817--1901), who retired the same
year. Friedrich Schur took his doctoral degree at the University
of Berlin with the thesis on the geometry of second order line
complexes. Later Friedrich Schur studied the algebraic surfaces of
third and fourth order and he made a valuable contribution to the
development of differential geometry by stating the famous theorem
on the spaces of constant curvature now bearing his name.
Friedrich Schur spent only four years (1888--1892) at the
University of Tartu. In this period he began to study the groups
of continuous transformations, then the rapidly developing new
branch of the differential geometry and due to his achievements in
this field, he can be reckoned a founder of this area of geometry
among such mathematicians as S. Lie, F. Engel, L. Maurer and W.
Killing.

Schur also studied the foundations of geometry. This is another
trend of research which played an essential role in the
development of geometry in Estonia. We only mention few names of
the geometers of the University of Tartu whose scientific activity
was related with the investigations on the foundations of
geometry.

Leonid Lachtin (1863--1927) was a professor of pure mathematics at
the University of Tartu for three years (1892--1895) succeeding
Friedrich Schur to the post. He studied the Lobachevskian
geometry and during his stay in Tartu published two papers
devoted to this subject, studying in one of them the Poincar{\'e}
interpretation of the Lobachevskian geometry. The successor of
Lachtin to the post of professor of pure mathematics V. Alekseyev
(1866--1943) studied the theory of the line congruences in
connection with surface theory and the theory of rational
invariants of bilinear forms.

The next 20th century is very significant in the history of
Estonia because it was marked by the independence of Estonia. We
touch gently the period from 1919 to 1950 because the special
stress in this period was laid on the investigations of the
foundations of geometry. The first professorship of mathematics at
the University of Tartu in the independent Estonia was occupied by
the Finnish mathematician Kalle V\"ais\"al\"a (1893--1968). He
spent in Tartu three years (1919--1922) and then moved to Turku
(Finland), where he took up the post of the professor of
mathematics at the University of Turku. In the period from 1930 to
1940 the investigations of the foundations of geometry at the
University of Tartu were continued in the papers of Jaan Sarv
(1877--1954), Arnold Humal (1908--1987), J\"uri Nuut (1892--1952) and
it should be mentioned that the approach developed by these
geometers was based on the notion of "betweenness", which is a
ternary relation on the set of points of a straight line
expressing the fact that one point lies between two others.

\section{Minimal surfaces and semiparallel submanifolds}

In this section we proceed to the next period of the development
of differential geometry in Estonia. This period covers the space
of time from 1950s to the beginning of 1990s, and the development
of differential geometry in this period in a great degree was
influenced by \"Ulo Lumiste. Let us mark the main stages of \"Ulo
Lumiste's scientific activity \cite{Lumiste1}. Lumiste was born in
V\"andra (Estonia) in 1929. He graduated from the University of Tartu
in 1952 and then he was sent to Moscow for post-graduate studies
at the Moscow University. In Moscow under the scientific
supervision of professor A. Vassiliev, Lumiste completed his thesis
devoted to the study of the geometry of submanifolds with fields
of asymptotic multidimensional directions in space forms and,
having successfully defended it, he obtained the candidate degree
(equivalent to PhD). Then he was appointed to a post of assistant
professor at the department of geometry of the University of
Tartu. In 1963--1965 Lumiste held a position of post-doctoral
researcher at the Moscow University. During this time he attended
the seminars on differential geometry led by S.P. Finikov, G.F.
Laptev and A.M. Vassiliev. Under the influence of these seminars
he began to study the theory of connections in fibre bundles and
its applications to geometry of families of homogeneous subspaces.
These investigations formed the basis for his DSc thesis, and he
defended it at the University of Kazan in 1968. In 1969 Lumiste
was appointed to the post of professor at the department of
algebra and geometry of the University of Tartu. He retired in
1995 and at present Lumiste is a professor emeritus.

Lumiste initiated a research in the
following areas of differential geometry:%
\begin{enumerate}
\item{the minimal and ruled surfaces, their generalizations;}
\item{canonical fibre bundles, induced connections and general
theory of connections;}
\item{2-parallel and semiparallel submanifolds;}
\item{connections in gauge theories.}
\end{enumerate}
Lumiste considerably contributed to each field of research from
mentioned above and he also continued the investigations on the
foundations of geometry started by his predecessors.

It is not possible to describe in a full extent the achievements
of \"U. Lumiste within the limits of this paper, therefore we
shall give only a brief description of the contribution of Lumiste
to some fields mentioned above and let us begin with the theory of
minimal surfaces. Lumiste showed in \cite{Lumiste2} that a minimal
surface of constant Gaussian curvature (other than plane) exists
only in the case of an elliptic space $S_n(c)$. He found all such
surfaces for $n\leq 5$ and proved that the Gaussian curvature of a
surface of this kind is $\frac{1c}{3}$ in the case of $S_4(c)$
(this is a so called Veronese surface) and the curvature of a
minimal surface with constant Gaussian curvature equals to zero in
the case of the elliptic spaces $S_3(c)$ and $S_5(c)$. It was also
shown in the same paper that every minimal surface of constant
Gaussian curvature is an orbit of a Lie group of isometries of a
corresponding elliptic space. In \cite{Lumiste3} Lumiste proved
the fundamental theorem for minimal surfaces, which states that a
minimal surface is entirely determined by its inner metric, the
principal curvatures and the angles between the principal
directions of all orders. It was also shown in the same paper that
every minimal surface can be bent continuously  within its own
class by leaving the values of first order principal curvatures
fixed. In \cite{Lumiste4} Lumiste proved that each indicatrix of
normal curvature of order $l$ of a minimal surface is a circle if and only if
the submanifold generated by the osculating planes of order l--1
is minimal.

This direction of research turned out to be very fruitful and
Lumiste drew his student L. Tuulmets (b. 1933) in the
investigations of the classes of 3-dimensional ruled surfaces in
the 4-dimensional space ${R}^4$, where ${R}^4$ can be either the
Euclidean space ${E}^4$ or the Minkowski space ${R}^{1,3}. $ The
starting point for Tuulmets' investigations was the paper
\cite{Lumiste4a}, where Lumiste elaborated the general theory of
quasi-congruences in the Euclidean space ${E}^4$. Tuulmets also
studied the various classes of surfaces in the 4-dimensional
Euclidean space ${E}^4$, Minkowski space ${R}^{1,3}$, projective
spaces $P^n$ and the spaces of constant curvature, where she used
the method of Cartan exterior forms and the systems of Pfaff
equations \cite{Tuulmets}. The question of consistency of a system
of Pfaff equations is very crucial in this kind of investigations,
but it also requires a large volume of pure algebraic computations
and Tuulmets assisted by an expert in the algebraic computer
methods R. Roomeldi successfully applied computer methods in her
investigations of Pfaff systems.

The congruences of null straight lines in the Minkowski space
$R^{1,3}$ were studied by R. Kolde (b. 1938) \cite{Kolde}. He
defined elliptic, hyperbolic and parabolic congruences of null
straight lines with the help of the local properties of
congruences. Using the notion of a central hypersurface, Kolde
constructed the set of symmetric tensors related to the rays in
the second order differential neighbourhoods. He showed that these
tensors could be used to canonize the frames in the case of the
congruences of hyperbolic and elliptic types. In the case of a
normal congruences it appeared that there were two different
canonical frames. In the special case of elliptic congruences
which are called isotropic similar canonical frames are determined
up to a parameter. Kolde found the geometrical meaning of these
canonical frames.

The next field of research initiated by \"U. Lumiste is the theory
of canonical fibre bundles, induced connections and the general
theory of connections \cite{Lumiste15} and \cite{Lumiste16}. A Grassmannian manifold of $m$-
dimensional
subspaces in Euclidean space or symplectic space can be considered
as a fiber bundle, whose fibers are $m$-dimensional subspaces.
There is a connection on this fiber bundle which can be defined in
a natural way. This connection allows to study the geometry of a
Grassmannian manifold by means of the curvature and the torsion of
the mentioned above connection. This direction of research was
developed by Lumiste's post-graduate students A. Parring (b.
1940), E. Abel (b. 1947) and A. Fleischer (b. 1948). Parring
studied a family of $2m$-dimensional symplectic subspaces in a
$2n$-dimensional affine-symplectic space interpreting this family
as a fiber bundle, whose standard fiber is a symplectic subspace
of this family \cite{Parring}. The group of symplectic motions
acts on each fiber of this fiber bundle. Parring used the
curvature and the torsion of an inner connection to study the
geometry of this fiber bundle and he also derived the structure
equations of this family of symplectic subspaces. Abel considered
a $(m+r)$-dimensional surface $V_{m+r}$ in the non-Euclidean space
\cite{Abel}. It can be shown that a surface $V_{m+r}$ stratifies,
where fibers are $r$-parametric families of $m$-dimensional
non-Euclidean subspaces. Elaborating the ideas presented in the
papers \cite{Lumiste15}, \cite{Parring} Abel defined three
connections on a surface $V_{m+r}$ and studied the properties of
the torsion and the curvature of these connections. Fleischer
studied homogeneous quotient spaces of the group of motions in
Euclidean space $E^4$ and the Lie triple systems
\cite{Fleischer2}. He also studied relations between reducibility
of the holonomy group $\text{hol}(\nabla)$ and properties of the
nonassociative algebra $m$ with multiplication defined by the
tensor $A$ \cite{Fleischer1}. Fleischer proved that if $M = G/H$
is a Riemannian non-symmetric reductive homogeneous space of a
simple Lie group $G$ and the isotropy representation of $H$ has
only inequivalent irreducible components, then any invariant
metric connection on $M$ has irreducible holonomy group.

Now we go on to the next field of research initiated by Lumiste
which is the theory of semiparallel submanifolds. Given a
$m$-dimensional smooth manifold $M^m$ immersed isometrically into
a $n$-dimensional Euclidean space $E^n$ one has two quadratic
differential forms associated with this submanifold, where one of
them is the first fundamental form $g=g_{ij}\,dx^{i}dx^{j},\;
i,j=1,2,\ldots,m$, where $x^1, x^2,\ldots, x^m$ are the local
coordinates of $M^m$, determined by the metric $(g_{ij})$, and the
other is the second fundamental form $h=h_{ij}\,dx^{i}dx^{j}$ with
values in the normal vector bundle over a submanifold $M^m$ (the
fiber of this vector bundle at a point $p$ of a submanifold $M^m$
is the orthogonal complement of the tangent space $T_p M^m$ of a
submanifold $M^m$ with respect to the whole Euclidean space
$E^n$). It is well known that the first fundamental form $g$ is
always parallel, i.e., $\nabla g=0$, where $\nabla$ is the
Levi-Civita connection, but the second fundamental form $h$ does
not need to be parallel. We remind that a submanifold $M^m$ in a
Euclidean space $E^n$ is said to be a parallel submanifold if
${\bar\nabla} h=0$, where ${\bar\nabla}$ is the van der
Waerden-Bortolotti connection on a submanifold $M^m$ which is a
pair of the Levi-Civita connection $\nabla$ and the normal
connection $\nabla^{\bot}$, i.e.,
$\bar\nabla=\nabla\oplus\nabla^{\bot}$. If
$\{e_{\alpha}\}$, where $\alpha=m+1,\ldots,n$, is an adapted orthonormal
local frame of the normal vector bundle, then
$h_{ij}=h_{ij}^{\alpha}\,e_{\alpha}$ and the components
$h_{ijk}^{\alpha}$ of ${\bar\nabla} h$ determined by
$h^{\alpha}_{ijk}= {\bar\nabla}_i h^{\alpha}_{kj}$ can be
expressed in terms of the components of the second fundamental
(mixed) tensor $h_{ij}^{\alpha}$ as follows:
\begin{equation}
h_{ijk}^{\alpha}\,\omega^k=%
    dh_{ij}^{\alpha}-h_{kj}^{\alpha}\,\omega_i^k-%
      h_{ik}^{\alpha}\,\omega_j^k+h_{ij}^{\beta}\,\omega_{\beta}^{\alpha},
\label{second}
\end{equation}
where $\omega^j_i$ are the connection 1-forms of the Levi-Civita
connection  $\nabla$ and $\omega^{\alpha}_{\beta}$ the connection
1-forms of the normal connection $\nabla^{\bot}$. In the case of
a parallel submanifold $M^m$ the components of ${\bar\nabla} h$
are all equal to zero, i.e., $h_{ijk}^{\alpha}=0$. Applying the
exterior differential to the both sides of (\ref{second}), we
obtain
\begin{equation}
{\bar\nabla} h_{ijk}^{\alpha}\wedge \omega^k=%
       h_{ij}^{\beta}\,\Omega^{\alpha}_{\beta}-%
         h_{kj}^{\alpha}\,\Omega^{k}_{i}-%
           h_{ik}^{\alpha}\,\Omega^{k}_{j},
\label{integrability1}
\end{equation}
where $\Omega^{k}_{i}$ is the curvature 2-form of the Levi-Civita
connection $\nabla$ and $\Omega^{\alpha}_{\beta}$ is the curvature
2-form of the normal connection $\nabla^{\bot}$. The above
relation (\ref{integrability1}) gives us the integrability
condition for the differential system ${\bar\nabla} h=0$ which is
\begin{equation}
h_{ij}^{\beta}\,\Omega^{\alpha}_{\beta}-%
         h_{kj}^{\alpha}\,\Omega^{k}_{i}-%
           h_{ik}^{\alpha}\,\Omega^{k}_{j}=0.%
\label{integrability2}
\end{equation}
A submanifold $M^m$ is said to be a semiparallel submanifold if
the integrability condition (\ref{integrability2}) is satisfied.

The term {\it semiparallel submanifolds} for submanifolds with
second fundamental form satisfying (\ref{integrability2}) was
introduced by J. Deprez in 1985. At the same time Lumiste together
with his post-graduate student V. Mirzoyan independently began to
study an interesting and important class of submanifolds with
parallel $\bar{\nabla} h$, which form the subclass in the class of
the semiparallel submanifolds as it follows from
(\ref{integrability2}). It should be noted that Deprez only
classified and described the semiparallel surfaces $M_2$ and
hypersurfaces $M_{n-1}$ in $E_n$ while in the main the theory of
semiparallel surfaces was developed by Lumiste. The theorem
asserting that any semiparallel submanifold $M_m$ in a space form
$M_n(c)$ is the second order envelope of an orbit of a Lie group
of isometries was proved by Lumiste in 1990 and this theorem plays
a crucial role in the theory of semiparallel submanifolds since it
suggests that in order to develop the theory of semiparallel
submanifolds one should first of all to find all the symmetric
orbits and then to describe the second order envelopes of these
orbits. It can be shown that any symmetric orbit is an orthogonal
product of irreducible ones, where an irreducible orbit is a
minimal submanifold in its sphere (except the case when
irreducible orbit is a plane). Lumiste showed that some
irreducible orbits can be constructed by means of mappings which
are known in the algebraic geometry. It turned out that symmetric
orbits of certain kind which arise in the connection with the
study of semiparallel submanifolds, whose first normal subspace at
any point has the maximal possible dimension $\frac{1}{2}m(m+1)$,
were earlier introduced and studied by R. Mullari (1931--1969)
\cite{Mullari1} and \cite{Mullari2} who was Lumiste's first
post-graduate student. Mullari considered a symmetric orbit which
is a space of constant curvature immersed into $E_n$ in such a way
that all its inner motions are generated by isometries of $E_n$
and $n=\frac{1}{2}m(m+3)$, where $m$ is the dimension of an orbit.
He called this kind of symmetric orbit {\it maximal symmetric}.
Later they were given the name {\it Veronese orbits} because any
of these orbits can be constructed as the image of the
$m$-dimensional sphere $S_m$ with respect to Veronese mapping. The
theorem proved by Lumiste \cite{Lumiste5} states that if $m\geq 2$
and $n=\frac{1}{2}m(m+3)$, then a complete semiparallel submanifold
$M_m$ in $E_n$ with maximal possible dimension $\frac{1}{2}m(m+1)$
of the first normal subspace at any point is a single Veronese
orbit.

The study of three-dimensional semiparallel submanifolds in $E_6$
(see \cite{Lumiste6} and \cite{Lumiste7}) led to a symmetric orbit, which
is generated by one-parameter family of\break 2-spheres, whose
orthogonal trajectories are circles. The direct generalization of
this is a \textit{Segre submanifold} $S_{(m_1,m_2)}$, which can be
constructed by means of the Segre map known in algebraic geometry.
The second order envelopes for Segre orbits can be found in
\cite{Lumiste8}. This direction of research is also developed by
K. Riives (b. 1942). Riives proved \cite{Riives1} that a
semiparallel submanifold $M^4$ of a Euclidean space $E^n$, $n>9$
is a second order envelope of $V^2(r_1) \times S^1(r_2)\times
S^1(r_3)$ (where $V^2 (r_1)$ is a Veronese surface in $E^5)$. In
the same paper some geometrical properties of Veronese and
Clifford leaves are described in terms of a certain function.
Using the Cartan moving frame method and the exterior differential
calculus Riives described in \cite{Riives2} some special classes of
curves on irreducible envelopes of the reducible symmetric
submanifolds $V^2(r_1) \times S^1 (r_2) \times S^1 (r_3)$ with a
Veronese component, which is a Veronese surface in $E^5$.

The third class of symmetric orbits, which were found by Lumiste
in the connection with the study of semiparallel submanifolds and
can be constructed with the help of a mapping known in algebraic
geometry. Let $G_{k,l}$ be the Grassmannian\break ($k$-dimensional
subspaces of the real Euclidean space $E_l$) and $T^{(0,k)}\subset
(E_l^{*})^{\otimes k}$ be the space of skew-symmetric
$k$-covariant tensors. We can consider the Grassmannian $G_{k,l}$
as a submanifold in $T^{(0,k)}$, where it turns to be an orbit. It
should be noted that the geometry of $G_{k,l}$ was also studied by
Lumiste's post-graduate student I. Maasikas (b. 1944) in \cite{Maasikas}.
Lumiste proved in \cite{Lumiste9} that this orbit is a symmetric
orbit only in the case of $k=2$ (it is called the Pl\"ucker
orbit). Later he showed that the second order envelope of Pl\"ucker
orbit $G_{2,l}$ is trivial, which means that it is neither more
nor less than the same orbit. Lumiste introduced the term {\it
umbilic-like orbits} for the symmetric orbits, whose second order
envelopes are trivial. He has shown that the class of umbilic-like
orbits includes not only the previously mentioned Pl\"ucker
orbits, but also the unitary orbits and Veronese-Grassmann orbits.
These results are summarized by Lumiste in his monographs
\cite{Lumiste12} and \cite{Lumiste8}.

Lumiste's post-graduate student M. V\"aljas (b. 1958) studied totally quasiumbilical
submanifolds with
nonflat normal connection and he proved in joint paper with \"U.
Lumiste \cite{Valjas1} the existence of totally quasiumbilical
submanifolds in Euclidean spaces with codimension 2 having a
non-flat normal connection. He also studied Dupin-Mannheim
submanifolds and Clifford cones in $E\sb{n+m}$ \cite{Valjas2}.
The other Lumiste's post-graduate student T. Vorovere (b. 1957) studied the evolute of
order $\lambda$ for a submanifold in a Euclidean space. A
submanifold $N$ in $E\sp n$ is called the evolute of order
$\lambda$ for a submanifold $M\sp m$ in $E\sp n$, if there exists
a submersion $f:N\to M$ such that the osculating plate of order
$(\lambda-1)$ of $N$ at a point $y\in N$ belong to the normal
space to $M$ at $f(x)$, and moreover, the osculating plane
intersects $M$ at $f(x)$. Virovere derived several criteria for
the existence of evolutes (see \cite{Virovere}).

\section{Fiber bundles, connections and gauge theories}

The question about the nature of a space and time has always
stirred up the minds of philosophers, mathematicians and
physicists. The geometry provides us with various mathematical
models of a space and this is the reason why the development of
geometry has been so closely related to the development of
classical mechanics and theoretical physics. The appearance of
non-abelian gauge theories in 50s of the 20th century gave a
powerful impetus to this relation. It is well known
that the theory of connections on principal and vector bundles is
an adequate geometric framework for the non-abelian gauge
theories. There were favorable conditions in Tartu in 70s to
elaborate the relation between theory of connections and gauge
theories because there was a group of physicists at the Institute
of Physics studying the gauge theories and a group of geometers at
the Department of Geometry and Algebra of the University of Tartu
studying the theory of connections on fiber bundles. The
physicists felt the need of basic knowledge in differential
geometry of fiber bundles such as connection, curvature,
characteristic classes (or Chern classes) and this led to the
joint seminar with geometers, which ran for several years (usually
once a week). This seminar was initiated by Madis K\~oiv on the
part of physicists and by \"Ulo Lumiste on the part of geometers.
This collaboration turned out to be very fruitful and it had not
only educational importance but it also led to an original
research in the areas such as Backlund transformations, quantum
Yang-Mills theory, supersymmetry, supermanifolds and supergravity.
In this paper we shall briefly describe the investigations
concerning a geometrical meaning of the Faddeev-Popov ghost fields
and BRST-transformations.
The first non-abelian gauge field theory with the gauge group
$SU(2)$ was constructed by Yang and Mills (see \cite{Yang-Mills}).
From a geometric point of view an adequate space for Yang-Mills
theory is a principal fiber bundle $P(M,G)$, where $M$ is the
smooth four dimensional manifold and $G=SU(2)$ is the structure
group of $P$ called a gauge group in field theory. A connection
1-form $\omega$ determines the Yang-Mills field potentials. Indeed
let $\pi: P\to M$ be a projection of a principal fiber bundle $P$
and $\sigma: V\to \pi^{-1}(V)$ be a section of some local
trivialization of $P$ over an open subset $V$ of $M$. Then the
pull-back of a connection form $\sigma^{*}\omega$ is the Lie
algebra $su(2)$-valued 1-form on a subset $V$ and can be expressed
in local coordinates $x^{\mu}$ of $V$ as follows
$\sigma^{*}\omega=A^{(V)}_{\mu}\,dx^{\mu}$. If the coefficients
$A^{(V)}_{\mu}$ satisfy the Yang-Mills equations one can interpret
them as Yang-Mills field potentials.  A connection 1-form $\omega$
with coefficients of its pull-back $\sigma^{*}\omega$ satisfying
the Yang-Mills equations is called a {\it Yang-Mills connection}
and its pull-back $\sigma^{*}\omega$ a {\it Yang-Mills} 1-{\it form}. If
$U$ is some other open subset of $M$ such that $U\cap V\neq
\emptyset$ and the section $\xi:U\to \pi^{-1}(U)$ of a local
trivialization over $U$ is related to the section $\sigma$ by the
$G$-valued function $g(x)$, i.e.,
\begin{equation}
\xi (x)=\sigma(x)\,g(x),\qquad x\in U\cap V, \nonumber
\end{equation}
then the pull-back $\xi^{*}\omega$ of a connection on $U$ is
related to the pull-back $\sigma^{*}\omega$ as follows
\begin{equation}
\xi^{*}\omega=Ad(g^{-1}(x))\,\sigma^{*}\omega+
                  g^{-1}(x)\,dg(x).
\nonumber
\end{equation}
If $\xi^{*}\omega=A^{(U)}_{\mu}\,dx^{\mu}$, then the above relation
leads to the relation between the corresponding Yang-Mills
potentials
\begin{equation}
A^{(U)}_{\mu}=g^{-1}(x)\,A^{(V)}_{\mu}\,g(x)+g^{-1}(x)\,\partial_{\mu}
g(x), \label{nonabelian}
\end{equation}
which is the gauge transformation in a non-abelian case. The
curvature 2-form
\begin{equation}
\Omega=D_{\omega}\,\omega=
           d\,\omega+\frac{1}{2}\,\lbrack \omega,\omega \rbrack,
\label{curvature}
\end{equation}
where $D_{\omega}$ is the covariant differential, written in local
coordinates
\begin{equation}
\Omega=F_{\mu\nu}\,dx^{\mu}dx^{\nu}, \label{strength}
\end{equation}
determines the strength $F_{\mu\nu}$ of a Yang-Mills field, which
can be expressed in the terms of potentials as follows
\begin{equation}
F_{\mu\nu}=\partial_{\mu} A_{\nu}-\partial_{\nu} A_{\mu}+
                               \lbrack A_{\mu},A_{\nu} \rbrack.
\end{equation}
The Yang-Mills action\vskip -.5cm
\begin{equation}
S_{YM}=\frac{1}{4}\,\int_M \langle \Omega,*\Omega\rangle,%
\label{action}
\end{equation}
where $*$ is the (Hodge) star operator and $\langle\, ,\, \rangle$
is a Killing form on the Lie algebra $su(2)$, is invariant with
respect to gauge transformations.

An invariance of the action $S_{YM}$ with respect to gauge
transformations (\ref{nonabelian}) brings a distinguishing feature
in the Yang-Mills field theory or, in a more general context, in
any non-abelian gauge field theory. Indeed, it shows that the real
physical configuration of a gauge field theory is determined not
by a single set of field potentials $A=\{A_{\mu}(x)\}$, but by the
entire class of gauge equivalent sets of field potentials, where
two sets $A=\{A_{\mu}(x)\}$ and $A'=\{{A'}_{\mu}(x)\}$ are said to be
gauge equivalent if\vskip -.5cm
\begin{equation}
A'_{\mu}=g^{-1}(x)\,A_{\mu}\,g(x)+g^{-1}(x)\,\partial_{\mu} g(x).
\end{equation}
Thus a real physical configuration is determined up to a gauge
transformation.  The complications which arise in a quantization
of a gauge theory are caused just by this ambiguity. This problem
can be fixed by a parametrization of the space of real physical
configurations where each real physical configuration will be
uniquely determined by a single set of parameters. This can be
done by imposing an additional condition on Yang-Mills potentials
and this condition should pick out a single representative from
every gauge equivalent class. This condition is called {\it
gauge-fixing condition} or simply gauge in gauge field theories.

The quantum Yang-Mills field theory was constructed by L. Faddeev,
V. Popov, B. De Witt and they showed that the approach of R.
Feynman based on functional integral was most suitable scheme for
quantization of gauge field theories because the principle of
gauge invariance could be expressed in terms of this approach very
easily: one should integrate not over the space of all field
configurations, but only over the space of gauge-equivalent
classes of field configurations. However applying this procedure
to the Yang-Mills field theory one encounters a problem of the
non-local functional ${\mbox det} M$, where $M$ is the operator
$$M(\alpha) =\partial^{\mu}\partial_{\mu}
         \alpha - \partial^{\mu} \lbrack A_{\mu}, \alpha\rbrack,$$
which appears in the functional integral for $S$-matrix. This was
the reason why the integrand in the functional integral of
generating functional was not in the canonical form $\exp(i\times
action)$. Faddeev and Popov solved this problem by
introducing the additional anticommuting fields $c(x)$ and ${\bar
c}(x)$, which allow to put the determinant ${\mbox det} M$ into a
form of the integral over the infinite dimensional Grassmann
algebra $\mathcal{R}$\  (here $c(x)$ and ${\bar c}(x)$ are the generators of
this algebra) as follows
\begin{equation}
{\mbox det} M= \int exp \lbrace i \int
     {\bar c}^a(x)\, M^{ab}\, c^b(x)\, dx \} \prod_x d{\bar c}dc.
\end{equation}
The anticommuting fields $c(x)$ and ${\bar c}(x)$ were later given the
name of {\it Faddeev-Popov ghost fields}. It was showed by C.
Becchi, A. Rouet, R. Stora and independently by I. Tyutin that the
quantum effective action (this is the action one obtains by adding
to the classical Yang-Mills action the terms generated by the
Faddeev-Popov ghost fields) is invariant with respect to
transformations
\begin{eqnarray}
A^a_{\mu}(x)&\to & A^a_{\mu}(x)+\nabla_{\mu}\,c(x)\,\epsilon,
\label{I}\\
c^a (x)&\to &
c^a(x)-\frac{1}{2}\,t^{a}_{bd}\,c^b(x)\,c^d(x)\,\epsilon,
\label{II}\\
{\bar c}^a (x)&\to & {\bar c}^a(x)+\lbrack
\partial_{\mu}A^a_{\mu}(x)\rbrack \,\epsilon,\label{III}
\end{eqnarray}
where $t^{a}_{bd}$ are the structure constants of the Lie algebra
$su(2)$, $\nabla_{\mu}$ is the covariant derivative and
$\epsilon$ is a Grassmann parameter $\epsilon^2=0$ anticommuting
with the Faddeev-Popov fields and commuting with the Yang-Mills
potentials. The transformations (\ref{I})--(\ref{III}) are called
the {\it BRST transformations}. This kind of transformations is
known in the modern field theory under the name of {\it
supersymmetries}. The remarkable property of the BRST
transformations is that they are nilpotent. The BRST
transformations induce the BRST operator
\begin{equation}
\delta\,A^a_{\mu}(x)= \nabla_{\mu}\,c^a(x),\;
     \delta\,c^a(x)=-\frac{1}{2}\,t^{a}_{bd}\,c^b(x)\,c^d(x),\;
        \delta\,{\bar c}^a(x)=\partial_{\mu}A^a_{\mu}(x),
        \label{BRST-operator}
\end{equation}
definition of BRST operator is usually replaced by an auxiliary
field $b^a(x)$ and the last formula in (\ref{BRST-operator})
takes on the form $\delta\,{\bar c}^a(x)=b^a(x)$. Later the
anti-BRST operator ${\bar\delta}$ was added to BRST operator and
together they form the {\it BRST-algebra}
\begin{equation}
\delta^2={\bar\delta}^2=0\;\;
     \delta\,{\bar\delta}+{\bar\delta}\,\delta=0.
\end{equation}
The appearance of the Faddeev-Popov ghost fields in the quantum
Yang-Mills theory raised an interesting problem of elaborating a
geometric structure which could  allow to incorporate these
additional fields into known geometric framework of gauge fields
based on fiber bundle technique and to derive the
BRST transformations from known equations. The property of
anticommutativity of the ghost fields suggests an idea to
construct their geometric interpretation by means of differential
forms since the wedge product of differential 1-forms  is also
skew-symmetric. This idea seems to be very alluring if we look at
the BRST transformation of the ghost fields (\ref{II}) which is
very similar to the Cartan-Maurer equation and this suggests to
construct the BRST operator by means of exterior derivative. The
geometric interpretation of the ghost fields and BRST
transformations based on the mentioned above idea was proposed and
developed in the papers \cite{Thierry-Mieg1} and \cite{Thierry-Mieg2}.
Though the geometric interpretation of the
ghost fields and BRST operators in terms of differential forms and
exterior derivative is very attractive, it has a problem with that
part of BRST operator $\delta$ which is determined by the
transformation of the anti-ghost field (\ref{III}). Indeed, the
BRST operator $\delta$ transforms the anti-ghost field ${\bar
c}^a(x)$ into an auxiliary bosonic field $b^a(x)$ and therefore
this part of the BRST operator is not consistent with the
properties of the exterior derivative which raises the degree of a
form by 1. In \cite{Lumiste10} Lumiste improved the interpretation
proposed in \cite{Thierry-Mieg1} and \cite{Thierry-Mieg2} and
extended it framework to include the anti-ghost fields ${\bar
c}(x)$ and the anti-BRST operator ${\bar\delta}$. The geometric
construction, he proposed for interpretation of the ghost field
$c(x)$, was not so rigid as in \cite{Thierry-Mieg1} and
\cite{Thierry-Mieg2} and this allowed to keep the dependence of
$c(x)$ on a gauge-fixing condition. He also showed that BRST
transformations for the fields $A_{\mu}(x)$ and $c(x)$ could be derived
from the well known Laptev equations for a connection on a
principal fiber bundle $P(M,G)$. Lumiste elaborated a formalism of
$q$-vector fields considered as $(-q)$-forms for geometric
interpretation of the anti-ghost field ${\bar c}(x)$ and an
analogue of the exterior derivative, which can be used together
with the well known operator $\star^{-1}\,d\,\star$, where $\star$
is the Hodge operator, for geometric interpretation of the
anti-BRST operator $\bar\delta$. Let us briefly describe the
geometric construction elaborated in \cite{Lumiste10}. Let $z$ be
a point of a principal fiber bundle $\pi: P \to M$ with a
structure group $G$, $p=\pi(z)\in M$ be the projection of $z$ and
$S_z$ be a subspace of the tangent space $T_zP$ such that
$T_zP=S_z\oplus V_zP$, where $V_zP$ is the tangent space to the
fiber $\pi^{-1}(p)$ passing through a point $z$. Let ${\mathcal
J}_P(z)$ be the set of all subspaces $S_z$ at a point $z$
satisfying $T_zP=S_z\oplus V_zP$. It can be proved that the set
${\mathcal J}_P=\bigcup_{z\in P} {\mathcal J}_P(z)$ is the smooth
manifold, which can be endowed with a structure of the principal
fiber bundle over $P$ with the projection $\pi': {\mathcal J}_P\to
P$ defined by $\pi'(S_z)=z$. If $R_g(z)=zg$ is the right action of
the group $G$ on a principal bundle $P$, then the right action
$R_g^{*}$ of $G$ on ${\mathcal J}_P$ is defined by
$R_g^{*}(S_z)=S_{zg}=dR_g(S_z)$, where $dR_g: T_zP\to T_{zg}P$ is
the differential of $R_g$.%

Let ${\mathcal V}(P)$ be the Lie algebra of vertical (or
fundamental) vector fields on $P$. It is well known that this Lie
algebra is isomorphic to the Lie algebra ${\underline G}$ of $G$,
i.e, ${\mathcal V}(P)\simeq {\underline G}$. If $Y$ is a
fundamental vector field then let us denote by ${\xi}_Y$ the
corresponding element of the Lie algebra ${\underline G}$ which
induces $Y$. Let $\Sigma: P\to {\mathcal J}_P$ be a smooth section
of the principal bundle ${\mathcal J_P}$. Any section $\Sigma$
generates the ${\underline G}$-valued $1$-form $\theta$ on $P$
which is defined as follows: if $X$ is a vector field on $P$, then
$X$ can be written as the sum of two vector fields $Y$ and $Z$, where
$Z_p\in \Sigma(p)$ and $Y$ is the uniquely determined vertical
vector field, and%
\begin{equation}
\theta(X)=\xi_{Y}.
\end{equation}%
The ${\underline G}$-valued $1$-form $\theta$ induces the
${\underline G}$-valued $1$-form $\tilde\theta$ on the principal
bundle ${\mathcal J}_P$ and the value of this form on a tangent
vector $\tilde v$ to ${\mathcal J}_P$ at a point $S_z$ is
determined by the formula
\begin{equation}
{\tilde\theta}(\tilde v)=\theta_z(d\pi'(\tilde v)).
\end{equation}
In \cite{Lumiste10} Lumiste proposed to consider the ${\underline
G}$-valued $1$-form $\tilde\theta$ as a geometric interpretation
for the Faddeev-Popov ghost field $c$. A smooth section $\Sigma$
used in the construction of $\tilde\theta$ can be interpreted in
terms of a gauge theory as a gauge-fixing condition and this shows
the advantage of the approach proposed by Lumiste which allows to
keep the dependence of the ghost field on a gauge.

The geometric interpretation of ghost fields and
BRST-supersymmetries in terms of differential forms suffers from a
shortcoming, which does not allow to take into account all
properties of the ghost fields. It is important that ghost fields
$c^a(x)$ and ${\bar c}^b(x)$ are the generators of an infinite
dimensional Grassmann algebra $\mathcal{R}$ as it is mentioned
above. This means that ghost fields anticommute not only with
respect to superscripts $a$ and $b$, but also with respect to a point $x$
of a base manifold $M$, i.e.,\vskip -.5cm
\begin{equation}
c^a(x)\,c^b(y)=-c^b(y)\,c^a(x),%
\label{ghostrelations}
\end{equation}
where $x$ and  $y$ are points of a base manifold $M$. The commutation
relations (\ref{ghostrelations}) show that if we consider the
ghost fields as generators of an infinite dimensional Grassmann
algebra, then the product of the ghost fields $c^a(x)$ and $c^b(y)$ is
determined even if $x$ and $y$ are different point of a manifold $M$. It
is well known that one can multiply differential forms pointwise
and the product of two differential forms has no sense if they are
taken in different points of a manifold. In order to incorporate
the anticommutation relations (\ref{ghostrelations}) into a
geometric interpretation Lumiste and the author of this paper
constructed an infinite dimensional supermanifold
$\mathcal{A}_{\mathcal{R}}$ (see \cite{Lumiste11}). There are few
approaches to a notion of a supermanifold and one of them was
proposed by F. Berezin. Briefly it can be described as follows:
given an ordinary smooth manifold $M$ one constructs a
supermanifold by means of the theory of sheaves, where $M$ is
usually called an underlying or a base manifold. The underlying
manifold $\mathcal{A}$ for an infinite dimensional supermanifold
$\mathcal{A}_{\mathcal{R}}$ proposed in \cite{Lumiste11} was the
infinite dimensional manifold of all smooth connections of a given
principal fiber bundle $P(M,G)$, where the differential structure
was defined with the help of a Banach space structure, and the
sheaves determining the structure of a supermanifold were
constructed with the help of infinite dimensional Grassmann
algebra $\mathcal{R}$. Figuratevily speaking this supermanifold
was infinite dimensional with respect to both sectors, even and
odd. The underlying manifold of all smooth connections
$\mathcal{A}$ has a rich geometric structure. It is well known
that it is an infinite dimensional principal fiber bundle with
respect to the action of the infinite dimensional Lie group of
gauge transformations and a gauge-fixing condition can be
considered as a section of this bundle. Given a smooth function
determined on an ordinary manifold $M$ one can extend it to the
smooth function on the supermanifold constructed over $M$. This
procedure is called a Grassmann analytical continuation in the
theory of supermanifolds. The Yang-Mills action (\ref{action}) can
be considered as a function on the infinite dimensional manifold
of the classes of gauge-equivalent connections and  the infinite
dimensional supermanifold $\mathcal{A}_{\mathcal{R}}$ proposed in
\cite{Lumiste11} allowed to show that the quantum effective action
is a Grassmann analytical continuation of the Yang-Mills action
(\ref{action}) along the fibers of the corresponding principal
fiber bundle (determined by the action of the group of gauge
transformations). The BRST-supersymmetries were interpreted as a
family of vector fields on the infinite dimensional supermanifold
$\mathcal{A}_{\mathcal{R}}$.

\section{Jet bundles and symmetries of differential equations}
Every time when we solve a differential equation, study the
singularities of a mapping or compute the invariants of a Lie
group we (in some way or other) use the structure of a jet space
$J_{n,m}$, which means that we make use of the operator of total
differentiation, Cartan' forms and Lie vector fields. This
direction of research in the area of differential geometry is
developed by Maido Rahula who succeeded to \"Ulo Lumiste on the
post of professor of geometry at the University of Tartu in 1990
when Lumiste retired.

Let us list the main periods of Maido Rahula biography. Rahula was
born in J\"arvamaa (Estonia) in 1936. When he was 13 years old his
family was banished to Siberia during the Stalin' deportations. He
graduated the University of Tomsk in 1959 and then came back to
Estonia, where he had the luck to enter the post-graduate studies
at the University of Tartu under the supervision of Lumiste though
according to the system existing at that time the members of
deported families were forbidden to enter the universities located
in the european part of the Soviet Union. Rahula defended his PhD
thesis "On higher order differential geometry" at the University
of Tartu in 1964. He spent four years (1967--1971) in Algeria
teaching mathematics at the University of Algeria within the
framework of the program of cooperation of the Soviet Union with
developing countries. The next period (1972--1989) of his life was
bound up with the Institute of Technology in Odessa.

Let us consider the jet space $J_{n,m}$ in the simplest case of
$n=m=1$. Let%
\begin{equation}
t, u, u', u'', \ldots
\end{equation}
be the coordinates in this space. The operator of total
differentiation
\begin{equation}
D=\frac{\partial}{\partial t}+u'\,\frac{\partial}{\partial u}+
               u''\,\frac{\partial}{\partial u'}+\ldots,
\end{equation}
can be considered as a linear vector field. The flow (i.e., the
one-parameter group of diffeomorphisms generated by this vector
field) is determined by the exponential law\vskip -.5cm
\begin{equation}
U'=C\,U \;\;{\rm and}\ \  U_t=e^{Ct}\,U, \label{exponentiallaw1}
\end{equation}
\vskip -.1cm
\noindent where $U$ is the column with infinite number of entries
$u,u',u'',\ldots $ and $C$ is the infinite-dimensional matrix,
whose non-zero elements can be obtained from the elements of the
main diagonal of the infinite-dimensional unit matrix by moving
each of them to the right along the corresponding row to the next
position. It is obvious that the matrix $C^l$, where $l$ is an
integer, can be obtained from the unit matrix by repeating $k$
times the previously described procedure. The exponential of the
matrix $C\,t$, where $t$ is a parameter, multiplied from the right
by the column $U$ yields the column $U_t$ whose first element\vskip -.5cm
\begin{equation}
u_t=\sum_{k=0}^{\infty} u^{(k)}\; \frac{t^k}{k!},
\end{equation}
\vskip -.2cm
\noindent and the next elements are the derivatives of $u_t$, i.e., $u'_t,
u''_t, \ldots $. Thus the formula $U_t=e^{Ct}\,U$ describes the
motion of a point $(t,U_t)$ along the trajectory of the vector
field $D$. In a general case of a jet space $J_{n,m}$ we have the
system of operators $D_i,\,i=1,2,\ldots,n$ and the formula
$U_t=e^{Ct}\,U$ written by means of multi-indices determines the
$n$-dimensional orbits of the additive group ${\bf R}^n$.

Let us denote by $\frac{\partial}{\partial U}$ the matrix whose
single row consists of the elements
\begin{equation}
\frac{\partial}{\partial u},\,\frac{\partial}{\partial u'},\,
                  \frac{\partial}{\partial u''},\ldots
\end{equation}
Then the operator $D$ can be written in the form
\begin{equation}
D=\frac{\partial}{\partial t}+\frac{\partial}{\partial U}\,U'.
\end{equation}
The entries of the matrix $(\frac{\partial}{\partial t}
\frac{\partial}{\partial U})$ form the basis in the jet space
$J_{1,1}$ and the
entries\break\vskip -.47cm \noindent  of the matrix  $(\begin{array}{c} dt \\ dU\\
\end{array})$ form the dual basis. If we replace the first\break\vskip -.4cm \noindent
entry $\frac{\partial}{\partial
t}$ in the first matrix by the operator of total differentiation
$D$ then the first entry in the second matrix containing the
elements of dual basis should be replaced by the Cartan form
$\omega=dU-U'\,dt$. This follows from the formulae
\begin{equation}
(\begin{array}{cc}\!\! D & \frac{\partial}{\partial U}\!\! \\ \end{array})
=(\begin{array}{cc}\!\!\frac{\partial}{\partial t} &
\frac{\partial}{\partial U}\!\! \\
\end{array})\,\left(%
\begin{array}{cc}
  1 & 0 \\
  U' & E \\
\end{array}%
\right),
\left(%
\begin{array}{c}
  dt \\
  \omega \\
\end{array}%
\right)=\left(%
\begin{array}{cc}
  1 & 0 \\
  -U' & E \\
\end{array}%
\right)\,\left(%
\begin{array}{c}
  dt \\
  dU \\
\end{array}%
\right). \nonumber
\end{equation}
The basis and its dual basis depend on a point of the jet space
and if this point starts to move along the trajectory of the
vector field $D$ passing through this point, then both basises
change and this change or dependence on a parameter $t$ can be
described by the same exponential law (\ref{exponentiallaw1}).
Indeed, if\vskip -.5cm
\begin{equation*}
    \left(\frac{\partial}{\partial U}\right)'=-\frac{\partial}{\partial
    U}\,C\;\;{\rm and}\ \
    \omega'=C\,\omega,
\nonumber
\end{equation*}
where the stroke stands for the Lie derivative with respect to
$D$, then
\begin{equation*}
\left(\frac{\partial}{\partial
U}\right)_t=\frac{\partial}{\partial
U}\,e^{-Ct}\;\;{\rm and}\ \ \omega_t=e^{Ct}\,\omega.
\end{equation*}

The formula $I=e^{-Ct}\,U$ determines an infinite number of
invariants of the operator $D$. Indeed, we have
$I'=e^{-Ct}\,(U'-CU)=0$ and taking into account
$dI=e^{-Ct}\,\omega$, we conclude that the exponential $e^{-Ct}$ is
an integrating matrix for the system of forms $\omega$. The dual
formula
\begin{equation*}
\frac{\partial}{\partial I}=\frac{\partial}{\partial U}\,e^{Ct}
\end{equation*}
determines an infinite number of invariant operators. These
operators form the basis for Lie vector fields and for
infinitesimal symmetries of the operator $D$.

Any vector field $P$ can be written in the above defined basises
as follows
\begin{equation}
P=\frac{\partial}{\partial t}\,\xi+\frac{\partial}{\partial
U}\,\lambda=D\,\xi+\frac{\partial}{\partial
U}\,\mu=D\,\xi+\frac{\partial}{\partial I}\,\nu,
\end{equation}
where
$\lambda=P\,U,\,\mu=\omega(D),\,\nu=P\,I,\,\mu=\lambda-U'\,\xi,\ {\rm and}\ \nu=e^{-
Ct}\,\mu$.
It can be proved that the following conditions are equivalent to
each other:
\begin{enumerate}
    \item vector field $P$ is a Lie vector field;
    \item $\nu'=0$;
    \item $\mu'=C\,\mu$;
    \item $\lambda'=C\,\lambda+U'\,\xi'$.
\end{enumerate}
The second condition is the simplest one, and it shows that the
components of $\nu$ are the invariants of the field $D$. The third
condition shows that the entries of the column $\mu$ are
$f,f',f'',\ldots$, where $f$ is a generating function. For
instance the functions $1,t,t^2/2,\ldots$ are the generating
functions for the field $\frac{\partial}{\partial I}$ which is a
vertical vector field since $\xi=0$. The fourth condition is more
complicated since a generating function does not enter it
explicitly, but this condition determines the components of a Lie
vector field in the natural basis. It should be mentioned that the
fourth condition can be found in the classical Lie theory.

The above described general scheme was developed by Rahula and its
more detailed description can be found in the monograph
\cite{Rahula1}. This scheme can be applied in the theory of
differential equations. Let $F(t,U)=0$ be a differential equation,
which we shall write in the form $F=0$. This equation determines
the surface $A_0$ in the jet space $J_{1,1}$. The stratification
of singularities $A_0\supset A_1\supset A_2\supset\ldots$ is
defined on the surface $A_0$ by means of the flow of the vector
field $D$, where $A_n$ is determined by the system of equations
$F^{(k)}=0$, where $k=0,1,2,\ldots,n$. The solutions of $F=0$ can be
constructed by means of those trajectories of the operator $D$
which belong to $A_n$ for each integer $n$. The integral of an
equation $F=0$ is determined by a function which is constant on
the surface $A_0$. The symmetries of $F=0$ are the transformations
of the jet space $J_{1,1}$ such that they leave invariant the
stratification arising on the surface $A_0$. The symmetries of
$F=0$ considered as mappings transform a solution of $F=0$ into
another solution of a same equation and they do the same thing
with the integrals of $F=0$. Thus we can conclude that the flow of
a vector field reproduces the solutions of a differential
equation.

It is useful to consider a Lie vector field as an extended group
operator. Any differential invariant (for example, the curvature of
a curve or the curvature of a surface) is the invariant of a Lie
vector field. Determination of symmetries in some sense an inverse
problem for the Erlanger Program of F. Klein (1872). Indeed, the
main aim of the Erlanger Program is to find the set of all
properties of a set $S$ that remain invariant when the elements of
this set are subjected to the transformations of some Lie group of
transformations.

An advantage of this approach is that it is based on a jet space
whose structure is universal. Indeed, let us consider the triple
$(D,t,U)$ in the jet space $J_{1,1}$, where $D$ is the operator of
total differentiation, $t$ is the canonical parameter for $D$ (it
can be interpreted as a time) and $U$ is the set of coordinates of
a fiber. Given the triple $(X,s,F)$, where $X$ is a vector field
defined on a manifold $M$, $s$ is the canonical parameter for $X$
and $F$ is the system of functions $X\,f,\;X^2\,f,\;\ldots$
generated by a smooth function $f$ defined on $M$, we can always
relate this triple to the triple $(D,t,U)$ with the help of a
mapping $\varphi: M\to J_{1,1}$  satisfying
\begin{equation}
t\circ \varphi=0,\;\; U\circ \varphi=F\ \  {\rm and}\ \
      (D I)\circ\varphi=X\,(I\circ\varphi),
\end{equation}
where $I$ is an arbitrary function defined on the space $J_{1,1}$.
In this way we can carry over any invariant of the operator $D$
from the jet space $J_{1,1}$ onto a manifold $M$, where it will be
the invariant of a vector field $X$. Particularly the invariants
$I\circ \varphi=e^{Cs}\,F$ on a manifold $M$ correspond to the
invariants $I=e^{Ct}\,U$ defined on the jet space. The covariant
tensor fields including the differential forms can be carried over
from the jet space $J_{1,1}$ onto a manifold $M$, where the
differential operators such as the Monge-Ampere operator,
Laplacian, Hessian, curvature operator and so on can be
constructed by means of these differential forms.

Thus we can use the structures defined on a jet space $J_{n,m}$ to
get a necessary information about the operators on a manifold $M$,
their invariants and symmetries. The set of all triples $(X,s,F)$
can be viewed as a category and the triple $(D,t,U)$ is a finite
object of this category. A universal problem is to construct an
initial element or a finite element of this category, and the
structure of a jet space $J_{n,m}$, which can be used to solve this
problem, is universal just in this sense. The structures briefly
described in this section, their analysis and possible
applications have been in detail described in the monograph
\cite{Rahula2}.

The geometric structures arising in the theory of differential
equations were also studied by H. Kilp (b. 1942)
in \cite{Kilp2} and \cite{Kilp3}.

\section{Noncommutative geometry and generalization of
supersymmetry}%
In this section we describe a direction of research in the field
of differential geometry, which was initiated by R. Kerner
(University Paris VI) at the beginning of 1990s and later
developed in cooperation with the colleagues from Paris, Wroclaw
(Poland) and the author of this paper. This direction is related
to the noncommutative geometry. During the last decade a
spectacular development of noncommutative generalizations of
differential geometry and Lie group theory has been achieved. The
respective new chapters of mathematical physics are known under
the names of {\it noncommutative geometry}, {\it quantum groups}
and {\it quantum spaces}. This section is based on the paper
\cite{Abramov1}.

Let us consider an associative algebra $\mathcal G$ over complex
numbers with generators $\theta ^A, A=1,2,..,N$ satisfying the
{\it ternary} relations
\begin{equation}
\theta ^A\theta ^B\theta ^C=j\theta ^B\theta ^C\theta ^A,%
\label{ternary relations}
\end{equation}
where $j$ is a primitive cube root of unity. We suppose that the
$N^2$ products $\theta ^A\theta ^B$ are linearly independent
entities. The algebra $\mathcal G$ with ternary commutation
relations (\ref{ternary relations}) is called a {\it ternary
Grassmann algebra} \cite{Abramov1}, because the commutation
relations (\ref{ternary relations}) are very similar to the
commutation relations of a classical Grassmann algebra. Indeed, if
$\theta^{\alpha}$ with $ \alpha=1,2,\ldots,n$ are the generators of the
Grassmann algebra with $n$ generators, then they are subjected to
the well known relations
$\theta^{\alpha}\theta^{\beta}=(-1)\theta^{\beta}\theta^{\alpha}$
which can be interpreted as follows: each permutation of
generators in the binary product $\theta^{\alpha}\theta^{\beta}$
is accompanied by multiplication by $-1$ and $-1$ is considered as
a primitive square root of unity. Thus one can get a ternary
analogue of the Grassmann algebra replacing a binary product of
generators by a ternary product, a permutation of two objects by a
cyclic permutation of three objects and a primitive square root of
unity by a primitive cube root of unity. It is obvious that
ternary analogue of Grassmann algebra is based on a faithful
representation of the cyclic group ${\Bbb Z}_3$ by cube roots of
unity.

Let us briefly describe the structure of the algebra $\mathcal G$.
The immediate corollary is that any product of four or more
generators must vanish. Here is the proof
\[
(\theta ^A\theta ^B\theta ^C)\theta ^D=j\theta ^B(\theta ^C\theta
^A\theta ^D)=j^2(\theta ^B\theta ^A\theta ^D)\theta ^C=\theta
^A(\theta ^D\theta ^B\theta ^C)=j\theta ^A\theta ^B\theta ^C\theta
^D.
\]
Now, as $(1-j)\not =0$, one must have $\theta ^A\theta ^B\theta
^C\theta^D=0 $. The dimension of the ternary Grassmann algebra
$\mathcal G$ is ${ N(N+1)(N+2)/3}+1$. Any cube of a generator is
equal to zero, i.e., $(\theta^A)^3=0$, and the odd permutation of
factors in a product of three leads to an independent quantity.

Our algebra admits a natural ${\Bbb Z}_3$-grading: under
multiplication, the grades add up modulo 3; the numbers are grade
0, the generators $\theta ^A$ are grade 1; the binary products are
grade 2 and the ternary products grade 0 again. The dimensions of
the subspaces of grade 0, 1 and 2 are, respectively, $N$ for grade
1, $N^2$ for grade 2 and $(N^3-N)/3+1$ for grade 0.

The lack of symmetry between the grades 1 and 2 (corresponding to
the generator $j$ and its square $j^2$ in the cyclic group ${\Bbb
Z}_3$, which are interchangeable, suggests that one should
introduce another set of $N$ generators of grade 2, whose squares
would be of grade 1, and which should obey the conjugate ternary
relations as follows
\[
\bar \theta ^{\bar A}\bar \theta ^{\bar B}\bar \theta ^{\bar
C}=j^2\bar \theta ^{\bar B}\bar \theta ^{\bar C}\bar \theta ^{\bar
A}.
\]
With respect to the ordinary generators $\theta ^A$, the conjugate
ones should behave like the products of two $\theta $'s, i.e.,
\begin{equation}
\theta ^A(\theta ^B\theta ^C)=j(\theta ^B\theta ^C)\theta
^A\rightarrow \theta ^A\bar \theta ^{\bar B}=j\bar \theta ^{\bar
B}\theta ^A \label{z2grad0}
\end{equation}
and consequently
\begin{equation}
\bar \theta ^{\bar B}\theta ^A=j^2\theta ^A\bar \theta ^{\bar B}.
\label{z2grad0bis}
\end{equation}
One may also note that there is an alternative choice for the
commutation relation between the ordinary and conjugate
generators, that makes the conjugate generators different from the
binary products of ordinary generators
\begin{equation}
\theta ^A\bar \theta ^{\bar B}=-j\bar \theta ^{\bar B}\theta ^A\text{ and }%
\bar \theta ^{\bar B}\theta ^A=-j^2\theta ^A\bar \theta ^{\bar B},
\label{z2grad1}
\end{equation}
which are still compatible with the ternary relations introduced
above.

This could be interpreted in the following way. We have assumed
that the algebra's field is the field of complex numbers, but we
can imagine that it is possible to multiply an element of the
${\Bbb Z}_3$-graded Grassmann algebra by an element of a {\it
binary} Grassmann algebra. We assume that the binary elements
commute with the ternary ones, but anticommute as usual with each
other. The ${\Bbb Z}_3$-graded Grassmann elements of a given grade
still have no binary commutation relation. Then our new algebra
admits two gradings: the ${\Bbb Z}_2$-grading and the ${\Bbb
Z}_3$-grading. The elements of ${\Bbb Z}_2$-grade 0 and ${\Bbb
Z}_3$-grades 1 and 2 obey the
rules (\ref{z2grad0}) and (\ref{z2grad0bis}) whereas the elements of ${\Bbb Z%
}_2$-grade 1 and ${\Bbb Z}_3$-grades 1 and 2 obey the rules
(\ref{z2grad1}). If we think that these objects can help in
modelling of the quark fields,
then a quark variable would be of ${\Bbb Z}_2$-grade 1 and ${\Bbb Z}_3$%
-grade 1, and an antiquark variable of ${\Bbb Z}_2$-grade 1 and ${\Bbb Z}_3$%
-grade 2. Then the products of a quark and an antiquark would
have both grades zero, making it a boson. In the same way, the
products of three quark
or three antiquark fields would be of ${\Bbb Z}_3$-grade 0 and of ${\Bbb Z}%
_2 $-grade 1, that is, they would very much look like a fermionic
field.

Now, the $\bar \theta $'s generate their own Grassmann subalgebra
of the same dimension that the one generated by $\theta $'s;
besides, we shall have all the mixed products containing both
types of generators, but which can be
always ordered e.g., with $\theta ^A$'s in front and $\bar \theta ^{\bar B}$%
's in the rear, by virtue of commutation relations. The products
of $\theta ^A$'s alone or of $\bar \theta ^{\bar A}$'s alone span
two subalgebras of
dimension $N(N+1)(N+2)/3$ each; the mixed products span new sectors of the $%
{\Bbb Z}_3$-graded Grassmann algebra.

In the case of usual ${\Bbb Z}_2$-graded Grassmann algebras the
anti-commutation between the generators of the algebra and the
assumed associativity imply automatically the fact that {\it all}
grade $0$ elements {\it commute} with the rest of the algebra,
while {\it any two} elements of grade $1$ anti-commute.

{In the case of the ${\Bbb Z}_3$-graded generalization such an
extension of ternary and binary relations {\it does not follow
automatically} and must be explicitly imposed. If we decide to
extend the relations (\ref{z2grad0}), (\ref{z2grad0bis}) and
(\ref{z2grad1}) to {\it all} elements of the algebra having a
well-defined grade (i.e., the monomials in $\theta $'s and $\bar
\theta $'s), then many additional expressions must vanish, e.g.,}
\[
{\theta ^A\underbrace{{\theta ^B{\bar \theta }}^{\bar C}}=\underbrace{{%
\theta ^B{\bar \theta }}^{\bar C}}\theta ^A=\theta
^B\underbrace{{{\bar \theta }^{\bar C}\theta }^{{A}}}={\bar \theta
}^C\theta ^A\theta ^B=0},
\]
{because on the one side, $\theta ^B{\bar \theta }^{\bar C}$ and
${\bar \theta }^{\bar C}\theta ^A$ are of grade 0 and commute with
all other elements, and on the other side, commuting ${\bar \theta
}^C$ with $\theta ^A\theta ^B$ one gets twice the factor $j$,
which leads to the overall factor $j^2{\bar \theta }^C\theta
^A\theta ^B$. This produces a
contradiction which can be solved only by supposing that $\theta ^A\theta ^B{%
\bar \theta }^C=0$.}

{The resulting ${\Bbb Z}_3$-graded algebra contains only the
following combinations of generators:
\vskip -.5cm
\[
A_1=\{\theta ,{\bar \theta }{\bar \theta }\}\text{, }A_2=\{{\bar \theta }%
,\theta \theta \}\ \ {\rm and}\ \ A_0=\{{\bf 1},\theta {\bar \theta
},\theta \theta \theta ,{\bar \theta }{\bar \theta }{\bar \theta
}\}.
\]}
The dimension of the algebra is
\[
D(N)=1{}+2N+3N^2+\frac{2(N^3-N)}3=\frac{3+4N+9N^2+2N^3}3.
\]
The four summands $1$, $2N$, $3N^2$ and $\frac{2(N^3-N)}3$
correspond to the subspaces respectively spanned by the
combinations $\{{\Bbb C}\}$, $\{\theta ,\bar \theta \}$, $\{\theta
\theta ,\theta \bar \theta ,\bar \theta \bar \theta \}$ and
$\{\theta \theta \theta ,\bar \theta \bar \theta \bar \theta \}$.

{Let us note that the set of grade $0$ (which obviously forms a
subalgebra of the ${\Bbb Z}_3$-graded Grassmann algebra) contains
the products which could symbolize the only observable
combinations of {\it quark fields} in quantum chromodynamics based
on the $SU(3)$-symmetry.}

We can introduce the ${\Bbb Z}_3$-graded derivations of the ${\Bbb Z}_3$%
-graded Grassmann algebra by postulating the following set of
rules\label {grassderiv}:
\[
\partial _A({\bf 1})=0\text{,\ \ }\partial _A\theta ^B={\delta }_A^B\ \ {\rm and}\
\ \partial _A\bar \theta ^{\bar B}=0
\]
and similarly\vskip -.5cm

\[
\partial _{\bar A}({\bf 1})=0\text{,\ }\ \partial _{\bar B}\bar \theta
^{\bar C}={\delta }_{\bar B}^{\bar C}\ \  {\rm and}\ \ \partial _{\bar
B}\theta ^A=0.
\]
When acting on various binary and ternary products, the derivation
rules are the following:
\[
\partial _A(\theta ^B\theta ^C)={\delta }_A^B\theta ^C+j{\delta }_A^C\theta
^B\ {\rm and}\ \partial _A(\theta ^B\theta ^C\theta ^D)={\delta
}_A^B\theta ^C\theta ^D+j{\delta }_A^C\theta ^D\theta
^B+j^2{\delta }_A^D\theta ^B\theta ^C\!.
\]
Similarly, for the conjugate entities,
\[
\partial _{\bar A}(\bar \theta ^{\bar B}\bar \theta ^{\bar C})={\delta }%
_{\bar A}^{\bar B}\bar \theta ^{\bar C}+j^2{\delta }_{\bar
A}^{\bar C}\bar \theta ^{\bar B}\ {\rm and}\ \partial _{\bar
A}(\bar \theta ^{\bar B}\bar \theta ^{\bar C}\bar \theta ^{\bar
D})={\delta }_{\bar A}^{\bar B}\bar \theta ^{\bar C}\bar \theta
^{\bar D}+j^2{\delta }_{\bar A}^{\bar C}\bar \theta ^{\bar D}\bar
\theta ^{\bar B}+j{\delta }_{\bar A}^{\bar D}\bar \theta ^{\bar
B}\bar \theta ^{\bar C}\!\!.
\]
We emphasize the "twisted" Leibniz rule for the ternary products in the above formulae.

Finally, for mixed binary products like $\theta ^A\bar \theta
^{\bar B}$, the derivation rules are the following:\vskip -.5cm
\[
\partial _A(\theta ^B\bar \theta ^{\bar C})={\delta }_A^B\bar \theta ^{\bar
C}\ \  {\rm and}\ \ {\partial }_{\bar A}(\theta ^B\bar \theta ^{\bar C})=j{\delta }%
_{\bar A}^{\bar C}{\theta }^B.
\]
There is no need for rules of derivation of fourth-order
homogeneous expressions, because these vanish identically.

As the immediate consequence of these rules, we have the following
important identities:
\[
\partial _A\partial _B\partial _C=j\partial _B\partial _C\partial _A\text{
and }\partial _{\bar A}\partial _{\bar B}\partial _{\bar
C}=j^2\partial _{\bar B}\partial _{\bar C}\partial _{\bar A},
\]
while\vskip -.5cm
\[
\partial _A\partial _{\bar C}=j\partial _{\bar C}\partial _A\text{ and }%
\partial _{\bar C}\partial _A=j^2\partial _A\partial _{\bar C}.
\]
Hence we have the important consequence
\begin{equation}
\partial _A\partial _B\partial _C+\partial _B\partial _C\partial _A+\partial
_C\partial _A\partial _B=0.  \label{sum3deriv}
\end{equation}
The ${\Bbb Z}_3$-graded generalization of the Grassmanian and the ${\Bbb Z}%
_3 $-graded derivatives defined above can be used in order to produce a $%
{\Bbb Z}_3$-generalization of the supersymmetry generators acting
on the usual ${\Bbb Z}_2$-graded Grassmann algebra generated by
anticommuting fermionic variables ${\theta }^\alpha $ and $\bar
\theta ^{\dot \beta }$ satifying the relations
\[
\theta ^\alpha \theta ^\beta +\theta ^\beta \theta ^\alpha =0\text{,\ \ }%
\bar \theta ^{\dot \alpha }\theta ^\beta +\theta ^\beta \bar
\theta ^{\dot \alpha }=0\ \  {\rm and}\ \ \bar \theta ^{\dot \alpha
}\bar \theta ^{\dot \beta }+\bar \theta ^{\dot \beta }\bar \theta
^{\dot \alpha }=0.
\]
The ``anti-Leibniz'' rule of derivation
\[
\partial _\alpha (\theta ^\beta \theta ^\gamma )=\delta _\alpha ^\beta
\theta ^\gamma -\delta _\alpha ^\gamma \theta ^\beta
\]
and similarly for any two dotted indices or mixed indices, one
verifies easily that all such derivations do anticommute:
\[
\partial _\alpha \partial _\beta +\partial _\beta \partial _\alpha =0\text{%
,\ \ }\partial _{\dot \alpha }\partial _{\dot \beta }+\partial
_{\dot \beta }\partial _{\dot \alpha }=0\ \  {\rm and}\ \ \partial
_\alpha \partial _{\dot \beta }+\partial _{\dot \beta }\partial
_\alpha =0. \label{sum2deriv}
\]
These rules enable us to construct the generators of the supersymmetric (or\break $%
{\Bbb Z}_2$-graded) ``odd'' translations
\[
{\mathcal
 D}_\alpha =\partial _\alpha +{\sigma }_{\alpha \dot\beta }^k\bar \theta ^{\dot \beta
}\partial _k\ \  {\rm and}\ \
{\mathcal D}_{\dot \beta }=\partial _{\dot \beta }+{\sigma }_{\alpha \dot \beta }^m{\theta
}^\alpha \partial _m,
\]
where both dotted and un-dotted indices $\alpha ,\dot \beta $ take
the values 1 and 2, while the space-time indices $k,l$ and $m$ run from 0
to 3. The anti-commutators of these differential operators yield
the ordinary (``even'') space-time translations
\[
{\mathcal D}_\alpha {\mathcal D}_{\dot \beta }+{\mathcal D}_{\dot
\beta }{\mathcal D}_\alpha =2\
{\sigma }_{\alpha \dot \beta
}^k\partial _k,
\]
while\vskip -.6
cm
\[
{\mathcal D}_\alpha {\mathcal D}_\beta +{\mathcal D}_\beta {\mathcal D}_\alpha =0\ \ {\rm
and}\ \ %
{\mathcal D}_{\dot \alpha }{\mathcal D}_{\dot \beta }+{\mathcal D}_{\dot \beta }{\mathcal
D}%
_{\dot \alpha }=0.
\]
The ${\Bbb Z}_3$-graded generalization would amount to find a
``cubic root'' of a linear differential operator, making use of
equation (\ref{sum3deriv}).
We must have six kinds of generalized Grassmann variables $\theta ^A$, $%
\theta ^{\stackrel{\wedge }{A}}$ and $\theta ^{\stackrel{\vee }{A}}$
on the one
hand and $\bar \theta ^{\bar A}$, $\bar \theta ^{\stackrel{\wedge }{\bar A}}$%
and $\bar \theta ^{\stackrel{\vee }{\bar A}}$ on the other hand,
which is formally analogous to the ${\Bbb Z}_2$-graded case.
Instead of the Pauli matrices we should introduce the entities
endowed with three indices (``cubic matrices'') with which the
generators of the ${\Bbb Z}_3$-graded translations of grade 1 and
2 may be constructed as follows:\vskip -.5cm
\[
{\it D}_A=\partial _A+{\rho }_{A\stackrel{\wedge }{B}\stackrel{\vee }{C}%
}^m\theta ^{\stackrel{\wedge }{B}}\theta ^{\stackrel{\vee }{C}}{\nabla }%
_m+\omega _{A\bar A}^m\bar \theta ^{\bar A}{\nabla }_m\text{,\ }{\it D}%
_{\bar A}=\partial _{\bar A}+{\bar \rho }_{\bar A\stackrel{\wedge }{\bar B}%
\stackrel{\vee }{\bar C}}^m\bar \theta ^{\stackrel{\wedge }{\bar
B}}\bar \theta ^{\stackrel{\vee }{\bar C}}{\nabla }_m+\bar \omega
_{\bar AA}^m\theta ^A{\nabla }_m,
\]
\vskip -.5cm
\[
{\it D}_{\stackrel{\wedge }{B}}=\partial _{\stackrel{\wedge }{B}}+{\rho }_{A%
\stackrel{\wedge }{B}\stackrel{\vee }{C}}^m\theta ^A\theta
^{\stackrel{\vee
}{C}}{\nabla }_m+\omega _{\stackrel{\wedge }{B}\stackrel{\wedge }{\bar B}%
}^m\bar \theta ^{\stackrel{\wedge }{\bar B}}{\nabla }_m\text{,\ }{\it D}_{%
\stackrel{\wedge }{\bar B}}=\partial _{\stackrel{\wedge }{\bar
B}}+{\bar \rho }_{\bar A\stackrel{\wedge }{\bar B}\stackrel{\vee
}{\bar C}}^m\bar \theta ^{\bar A}\bar \theta ^{\stackrel{\vee
}{\bar C}}{\nabla }_m+\bar
\omega _{\stackrel{\wedge }{\bar B}\stackrel{\wedge }{B}}^m\theta ^{%
\stackrel{\wedge }{B}}{\nabla }_m
\]\noindent and\vskip -.6cm
\[
{\it D}_{\stackrel{\vee }{C}}=\partial _{\stackrel{\vee }{C}}+{\rho }_{A%
\stackrel{\wedge }{B}\stackrel{\vee }{C}}^m\theta ^A\theta ^{\stackrel{%
\wedge }{B}}{\nabla }_m+\omega _{\stackrel{\vee }{C}\stackrel{\vee }{\bar C}%
}^m\bar \theta ^{\stackrel{\vee }{\bar C}}{\nabla }_m\text{,\ }{\it D}_{%
\stackrel{\vee }{\bar C}}=\partial _{\stackrel{\vee }{\bar C}}+{\bar \rho }%
_{\bar A\stackrel{\wedge }{\bar B}\stackrel{\vee }{\bar C}}^m\bar
\theta
^{\bar A}\bar \theta ^{\stackrel{\wedge }{\bar B}}{\nabla }_m+\bar \omega _{%
\stackrel{\vee }{\bar C}\stackrel{\vee }{C}}^m\theta ^{\stackrel{\vee }{C}}{%
\nabla }_m.
\]
The nature of the indices needs not to be specified; the only
important thing to be assumed at this stage is that the
differential operators $\nabla _m$ do commute with the ${\Bbb
Z}_3$-graded differentiations $\partial _A$. It is also
interesting to consider the operators one gets when the $\nabla _m
$ are replaced with {\em supersymmetric} derivations (that
anticommute with the ${\Bbb Z}_3$-graded differentiations). But in
the simpler case described here, the following operators acting on
the ${\Bbb Z}_3$-graded generalized Grassmanian:
$$D_{ABC}^{III}={\it D}_A{\it D}_B{\it D}_C+{\it D}_B{\it D}_C{\it D}_A+%
{\it D}_C{\it D}_A{\it D}_B+{\it D}_C{\it D}_B{\it D}_A+{\it D}_B{\it D}_A%
{\it D}_C+{\it D}_A{\it D}_C{\it D}_B, $$
$$\bar D_{\bar A\bar B\bar C}^{III} ={\it D}_{\bar A}{\it D}_{\bar B}{\it D}%
_{\bar C}+{\it D}_{\bar B}{\it D}_{\bar C}{\it D}_{\bar A}+{\it D}_{\bar C}%
{\it D}_{\bar A}{\it D}_{\bar B}+{\it D}_{\bar C}{\it D}_{\bar B}{\it D}%
_{\bar A}+{\it D}_{\bar B}{\it D}_{\bar A}{\it D}_{\bar C}+{\it D}_{\bar A}%
{\it D}_{\bar C}{\it D}_{\bar B}$$\noindent and
$$D_{A\bar A}^{II} ={\it D}_A{\it D}_{\bar A}-j^2{\it D}_{\bar
A}{\it D}_A$$
represent {\it homogeneous} operators on the ${\Bbb Z}_3$-graded
Grassmann algebra, i.e., they map polynomials in $\theta $'s of a
given grade into polynomials of the same grade; the result can be
represented by a complex-valued matrix containing various
combinations of the differentiations $\nabla _m$ ; their eventual
symmetry properties will
depend on the assumed symmetry properties of the matrices $\rho _{ABC}$ and $%
\omega _{A\bar B}$.

Let us consider in more detail the case of dimension $3$ (the
simplest possible realization of the ${\Bbb Z}_3$-graded
Grassmannian and the derivations on it is of course the case with
{\it one} generator and its conjugate.

The dimension of the ${\Bbb Z}_3$-graded Grassmann algebra with three grade-$%
1$ genera-\break\vskip -.5cm\noindent
tors $\theta $, $\stackrel{\wedge }{\theta }$ and $\stackrel{\vee }{%
\theta }$ and three ``conjugate'' grade-$2$ generators $\bar \theta $, $%
\stackrel{\wedge }{\bar \theta }$ and $\stackrel{\vee }{\bar \theta }$ is $%
51 $; any linear operator, including the derivations $\partial _A$
and the multiplication by any combination of the generators, as
well as the operators ${\it D}_A$ and ${\it D}_{\bar A}$
introduced above, can be represented by means of $51\times 51$
complex-valued matrices. Unfortunately, the operators $D^{II}$ and
$D^{III}$ are neither diagonal nor diagonalizable. But if we apply
them to a scalar function $f$, we get\vskip -.5cm
\[
D_{1\bar 1}^{II}f=(\omega _{1\bar 1}^m+\bar \omega _{\bar
11}^m)\nabla _mf
\]
\vskip -.2cm\noindent and\vskip -.5cm
\[
D_{1\stackrel{\wedge }{1}\stackrel{\vee }{1}}^{III}f=-3j^2\rho _{1\stackrel{%
\wedge }{1}\stackrel{\vee }{1}}^m\nabla _mf\text{ as well as }\bar D_{\bar 1%
\stackrel{\wedge }{\bar 1}\stackrel{\vee }{\bar 1}}^{III}f=-3j\bar
\rho _{\bar 1\stackrel{\wedge }{\bar 1}\stackrel{\vee }{\bar
1}}^m\nabla _mf.
\]
The $\omega $ matrices are the only ones that remain in the
$D^{II}$ whereas the $\rho $ cubic matrices emerge from the
ternary combinations $D^{III}$. On the space of scalar functions,
our operators act simultaneously{\em \ }as {\em square }and {\em
cubic} roots of ordinary translations. Using extensions of these
objects, where ${\nabla }_m$ are replaced with the
supersymmetry generators, we have contructed a simple ${\Bbb
Z}_3$-graded noncommutative geometry model featuring three Higgs
fields. The lagrangian contains the potential term\vskip -.6cm
\[
V=3\left| \Phi _1+\Phi _2+\Phi _3+\Phi _1\Phi _2+\Phi _2\Phi
_3+\Phi _3\Phi _1+\Phi _1\Phi _2\Phi _3\right| ^2
\]
and implies multiple spontaneous symmetry breaking.

The ternary
generalization of Grassmann algebra described in this section can
be used to construct a generalization of exterior calculus with
exterior differential $d$ satisfying $d^N=0$ with $N>2$. This
direction of research has been developed by Abramov's
post-graduste student N. Bazunova (b. 1964) in \cite{bazunova}.

\vskip.4cm%
\noindent%
{\bf Acknowledgments}%
\vskip.3cm%
I am grateful to the organizers of the Estonian-Finnish Seminar on
the Development of Mathematics for an invitation to give a talk on
the development of differential geometry in Estonia. I would like
to express my gratitude to \"U. Lumiste who took the trouble to
read the manuscript of this paper and his suggestions have been
very valuable for me in preparation of the final version of this
paper. I am also grateful to my colleagues M. Rahula, A. Parring,
L. Tuulmets, K. Riives, E. Abel, H. Kilp for the valuable
explanations concerning their research and the biographic
materials, which they kindly placed at my disposal. I gratefully
acknowledge the financial support of the Estonian Science
Foundation under the grant No. 4515.
\bibliographystyle{amsplain}

\end{document}